# Non-null Slant Ruled Surfaces

**Mehmet Önder**
*Independent Researcher, Delibekirli Village, Tepe Street, No. 63, 31440 Kırıkhan, Hatay, Turkey.*
E-mail: mehmetonder197999@gmail.com

**Abstract**
In this study, we define some new types of non-null ruled surfaces called slant ruled surfaces in the Minkowski 3-space $E_1^3$. We introduce some characterizations for a non-null ruled surface to be a slant ruled surface in $E_1^3$. Moreover, we obtain some corollaries which give the relationships between a non-null slant ruled surface and its striction line in $E_1^3$.

## 1. Introduction

In the study of curve theory, the curves whose curvatures satisfy some special conditions have an important role. The well-known of such curves is general helix defined by the classical definition that the tangent lines of the curve makes a constant angle with a fixed straight line [4]. In 1802, M.A. Lancret stated a result on the helices which was first proved by B. de Saint Venant in 1845 [20]. Venant showed that a curve is a general helix if and only if the ratio of the curvatures $\kappa$ and $\tau$ of the curve is constant, i.e., $\kappa/\tau$ is constant at all points of the curve. Helices have been studied not only in Euclidean spaces but also in Lorentzian spaces by some mathematicians and different characterizations of these curves have been obtained according to the properties of the spaces [6,7,10,15].

Recently, Izumiya and Takeuchi have introduced a new curve called slant helix which is defined by the property the normal lines of the curve make a constant angle with a fixed direction in Euclidean 3-space $E^3$ [8]. Later, the spherical images, the tangent indicatrix and the binormal indicatrix of a slant helix have been studied by Kula and Yaylı and they have obtained that the spherical images of a slant helix are spherical helices [11]. The position vector of a slant helix in $E^3$ has been studied by Ali [3]. Then the corresponding characterizations for the position vector of a timelike slant helix in Minkowski 3-space $E_1^3$ have been given by Ali and Turgut [2]. Recently, Ali and Lopez have also given some new characterizations of slant helices in Minkowski 3-space $E_1^3$ [1].

Analogue to the curves, ruled surfaces have orhonormal frames along their striction curves. So, the notion "helix" or "slant helix" can be considered for ruled surfaces. Before, Önder and Kaya have studied this subject for null scrolls and defined slant null scrolls in $E_1^3$ [19]. In this paper, we define non-null slant ruled surfaces by considering the Frenet vectors of timelike and spacelike ruled surfaces in $E_1^3$. We give the conditions for a non-null ruled surface to be a slant ruled surface.

## 2. Preliminaries

Let $E_1^3$ be a Minkowski 3-space with natural Lorentz Metric $\langle , \rangle = -dx_1^2 + dx_2^2 + dx_3^2$, where $(x_1, x_2, x_3)$ is a rectangular coordinate system of $E_1^3$. Since this metric is not positive definite, for an arbitrary vector $\vec{v} = (v_1, v_2, v_3)$ in $E_1^3$ we have i) $\langle \vec{v}, \vec{v} \rangle > 0$ and $\vec{v} = 0$, ii) $\langle \vec{v}, \vec{v} \rangle < 0$ iii) $\langle \vec{v}, \vec{v} \rangle = 0$ and $\vec{v} \neq 0$ [13]. Then we have three types of vectors: spacelike, timelike or null(lightlike) if (i), (ii) or (iii) holds, respectively. Similarly, an arbitrary curve $\vec{\alpha} = \vec{\alpha}(s)$ can locally be spacelike, timelike or null (lightlike), if all of its velocity vectors $\vec{\alpha}'(s)$ satisfy (i), (ii) or (iii), respectively. For the vectors $\vec{x} = (x_1, x_2, x_3)$ and $\vec{y} = (y_1, y_2, y_3)$ in $E_1^3$, the vector product of $\vec{x}$ and $\vec{y}$ is defined by

$$\vec{x} \times \vec{y} = (x_2 y_3 - x_3 y_2, x_1 y_3 - x_3 y_1, x_2 y_1 - x_1 y_2).$$

The Lorentzian sphere and hyperbolic sphere of radius $r$ and center origin 0 in $E_1^3$ are given by

$$S_1^2 = \{\vec{x} = (x_1, x_2, x_3) \in E_1^3 : \langle \vec{x}, \vec{x} \rangle = r^2\},$$

and

$$H_0^2 = \{\vec{x} = (x_1, x_2, x_3) \in E_1^3 : \langle \vec{x}, \vec{x} \rangle = -r^2\},$$

respectively [21].

Analogue to the curves, a surface can be timelike or spacelike in $E_1^3$. The Lorentzian character of a surface in $E_1^3$ is determined by the induced metric on the surface. The surface is called timelike(spacelike), if this metric is a Lorentz metric(positive definite Riemannian metric) [5,21].

Let now $I$ be an open interval in the real line $IR$. Let $\vec{k} = \vec{k}(u)$ be a curve in $E_1^3$ defined on $I$ and $\vec{q} = \vec{q}(u)$ be a unit direction vector of an oriented line in $E_1^3$. Then we have the following parametrization for a ruled surface $N$,

$$\vec{r}(u,v) = \vec{k}(u) + v\vec{q}(u). \tag{1}$$

The straight lines of the surface are called rulings and the curve $\vec{k} = \vec{k}(u)$ is called base curve or generating curve. In particular, if the direction of $\vec{q}$ is constant, the ruled surface is said to be cylindrical, and non-cylindrical otherwise.

The function defined by

$$\delta = \frac{|d\vec{k}, \vec{q}, d\vec{q}|}{\langle d\vec{q}, d\vec{q} \rangle}$$

is called the distribution parameter (or drall) of the ruled surface. Then, $N$ is called developable surface if and only if $\delta = 0$ [12,16,18]. Then at all points of same ruling, the tangent planes are identical, i.e., tangent plane contacts the surface along a ruling. If $|d\vec{k}, \vec{q}, d\vec{q}| \neq 0$, then the tangent planes of the surface $N$ are distinct at all points of same ruling which is called nontorsal [16,18].

Let consider the unit normal vector $\vec{m}$ of $N$ defined by $\vec{m} = \frac{\vec{r}_u \times \vec{r}_v}{\|\vec{r}_u \times \vec{r}_v\|}$. So, at the points of a nontorsal ruling $u = u_1$ we have

$$\vec{a} = \lim_{v \to \infty} \vec{m}(u_1, v) = \frac{d\vec{q} \times \vec{q}}{\|d\vec{q}\|}.$$

which is called central tangent. The point at which the vectors $\vec{a}$ and $\vec{m}$ are ortogonal is called the striction point (or central point) $C$ and the set of striction points of all rulings is called striction curve which has the parametric representation

$$\vec{c}(u) = \vec{k}(u) - \frac{\langle d\vec{q}, d\vec{k} \rangle}{\langle d\vec{q}, d\vec{q} \rangle} \vec{q}(u). \tag{2}$$

It is clear that the base curve is the same with striction curve if and only if $\langle d\vec{q}, d\vec{k} \rangle = 0$.

Since the vectors $\vec{a}$ and $\vec{q}$ are orthogonal, we can define an orthonormal frame on the surface. For this purpose, let write $\vec{h} = \vec{a} \times \vec{q}$. The unit vector $\vec{h}$ is called central normal and the orthonormal frame $\{C; \vec{q}, \vec{h}, \vec{a}\}$ at central point $C$ is called Frenet frame of $N$.

According to the Lorentzian casual characters of ruling and central normal, the Lorentzian character of the surface $N$ is classified as follows;

**i)** If the central normal vector $\vec{h}$ is spacelike and $\vec{q}$ is timelike, then the ruled surface $N$ is said to be of type $N_-$.

**ii)** If the central normal vector $\vec{h}$ and the ruling $\vec{q}$ are both spacelike, then the ruled surface $N$ is said to be of type $N_+$.

**iii)** If the central normal vector $\vec{h}$ is timelike, then the ruled surface $N$ is said to be of type $N_\times$ [16,18].

The ruled surfaces of type $N_-$ and $N_+$ are clearly timelike and the ruled surface of type $N_\times$ is spacelike. By using these classifications and taking the striction curve as the base curve the parametrization of the ruled surface $N$ can be given as follows,

$$\vec{r}(s,v) = \vec{c}(s) + v\vec{q}(s), \tag{3}$$

where $\langle \vec{q}, \vec{q} \rangle = \varepsilon (=\pm 1)$, $\langle \vec{h}, \vec{h} \rangle = \pm 1$ and $s$ is the arc length of the striction curve.

For the derivatives of the vectors of Frenet frame $\{C; \vec{q}, \vec{h}, \vec{a}\}$ of ruled surface $N$ with respect to the arc length $s$ of striction curve we have the followings

**i)** If the ruled surface $N$ is a timelike ruled surface then we have

$$\begin{bmatrix} d\vec{q}/ds \\ d\vec{h}/ds \\ d\vec{a}/ds \end{bmatrix} = \begin{bmatrix} 0 & k_1 & 0 \\ -\varepsilon k_1 & 0 & k_2 \\ 0 & \varepsilon k_2 & 0 \end{bmatrix} \begin{bmatrix} \vec{q} \\ \vec{h} \\ \vec{a} \end{bmatrix}, \tag{4}$$

and $\vec{q} \times \vec{h} = \varepsilon \vec{a}$, $\vec{h} \times \vec{a} = -\varepsilon \vec{q}$, $\vec{a} \times \vec{q} = -\vec{h}$ [16].

**ii)** If the ruled surface $N$ is spacelike ruled surface then we have

$$\begin{bmatrix} d\vec{q}/ds \\ d\vec{h}/ds \\ d\vec{a}/ds \end{bmatrix} = \begin{bmatrix} 0 & k_1 & 0 \\ k_1 & 0 & k_2 \\ 0 & k_2 & 0 \end{bmatrix} \begin{bmatrix} \vec{q} \\ \vec{h} \\ \vec{a} \end{bmatrix}, \tag{5}$$

and $\vec{q} \times \vec{h} = -\vec{a}$, $\vec{h} \times \vec{a} = -\vec{q}$, $\vec{a} \times \vec{q} = \vec{h}$ [18].

In the equations (4) and (5), $k_1 = \frac{ds_1}{ds}$, $k_2 = \frac{ds_3}{ds}$ and $s_1$, $s_3$ are the arc lengths of the spherical curves generated by the unit vectors $\vec{q}$ and $\vec{a}$, respectively.

**Theorem 2.1** ([14]). *Let the striction curve $\vec{c} = \vec{c}(s)$ of ruled surface $N$ be a unit speed curve and have the same Lorentzian casual character with the ruling and let $\vec{c}(s)$ be the base curve* *of the surface. Then $N$ is developable if and only if the unit tangent of the striction curve is the same with the ruling along the curve.*

## 3. $q$-Slant Ruled Surfaces in $E_1^3$

In this section, we introduce the definition and characterizations of $q$-slant ruled surfaces in $E_1^3$. First, we give the following definition.

**Definition 3.1.** Let $N$ be a non-null ruled surface in $E_1^3$ given by the parametrization
$$\vec{r}(s,v) = \vec{c}(s) + v\vec{q}(s), \quad \|\vec{q}(s)\| = 1, \tag{6}$$
where $\vec{c}(s)$ is striction curve of $N$ and $s$ is arc length parameter of $\vec{c}(s)$. Let the Frenet frame and non-zero invariants of $N$ be $\{\vec{q}, \vec{h}, \vec{a}\}$ and $k_1, k_2$, respectively. Then, $N$ is called a $q$-slant ruled surface if the ruling $\vec{q}(s)$ makes a constant angle with a fixed non-null unit direction $\vec{u}$ in the space, i.e.,
$$\langle \vec{q}, \vec{u} \rangle = c_q = constant. \tag{7}$$

Then we give the following characterizations for $q$-slant ruled surfaces in $E_1^3$. Whenever we talk about $N$ we will mean that the surface has the properties as assumed in Definition 3.1.

**Theorem 3.1.** *The ruled surface $N$ is a $q$-slant ruled surface if and only if the function $k_1/k_2$ is constant and given by*
$$k_1/k_2 = \begin{cases} -\varepsilon c_a/c_q, & N \text{ is timelike} \\ -c_a/c_q, & N \text{ is spacelike} \end{cases} \tag{8}$$
*where $c_q = \langle \vec{q}, \vec{u} \rangle$, $c_a = \langle \vec{a}, \vec{u} \rangle$ are nen-zero constants.*

**Proof:** Let $N$ be a $q$-slant ruled surface in $E_1^3$. Then denoting by $\vec{u}$ the unit vector of fixed direction, $N$ satisfies
$$\langle \vec{q}, \vec{u} \rangle = c_q = constant. \tag{9}$$
Differentiating (9) with respect to $s$ gives $\langle \vec{h}, \vec{u} \rangle = 0$. Therefore, $\vec{u}$ lies on the plane spanned by the vectors $\vec{q}$ and $\vec{a}$, i.e.,
$$\vec{u} = c_q \vec{q} + c_a \vec{a}, \tag{10}$$
where $c_q$ and $c_a$ are real constants. By differentiating (10) with respect to $s$ it follows
$$0 = \begin{cases} (c_q k_1 + \varepsilon c_a k_2)\vec{h}; & N \text{ is timelike,} \\ (c_q k_1 + c_a k_2)\vec{h}; & N \text{ is spacelike,} \end{cases} \tag{11}$$
and then we have that the function
$$k_1/k_2 = \begin{cases} -\varepsilon c_a/c_q, & N \text{ is timelike} \\ -c_a/c_q, & N \text{ is spacelike} \end{cases}$$
is constant.

Conversely, given a non-null ruled surface $N$, the equation (8) is satisfied. We define
$$\vec{u} = c_q \vec{q} + c_a \vec{a}, \tag{12}$$

where $\langle \vec{q}, \vec{u} \rangle = c_q$, $\langle \vec{a}, \vec{u} \rangle = c_a$ are non-zero constants. Differentiating (12) and using (8) it follows $\vec{u}' = 0$, i.e., $\vec{u}$ is a constant vector. On the other hand $\langle \vec{q}, \vec{u} \rangle = c_q = constant$. Then $N$ is a $q$-slant ruled surface in $E_1^3$.

***Theorem 3.2.*** *Non-null ruled surface $N$ is a $q$-slant ruled surface if and only if* $\det(\vec{q}', \vec{q}'', \vec{q}''') = 0$ *holds.*

**Proof:** From the Frenet formulae in (4) and (5) we have

$$\det(\vec{q}', \vec{q}'', \vec{q}''') = \begin{cases} -\varepsilon k_1^3 k_2^2 \left(\dfrac{k_1}{k_2}\right)'; & N \text{ is timelike} \\ k_1^3 k_2^2 \left(\dfrac{k_1}{k_2}\right)'; & N \text{ is spacelike} \end{cases} \quad (13)$$

Let now $N$ be a $q$-slant ruled surface in $E_1^3$. By Theorem 3.1 we have $k_1/k_2$ is constant. Then from (13) it follows that $\det(\vec{q}', \vec{q}'', \vec{q}''') = 0$.

Conversely, if $\det(\vec{q}', \vec{q}'', \vec{q}''') = 0$, since the curvatures are non-zero from (13) it is obtained that $k_1/k_2$ is constant and Theorem 3.1 gives that $N$ is a $q$-slant ruled surface in $E_1^3$.

***Theorem 3.3.*** *Non-null ruled surface $N$ is a $q$-slant ruled surface if and only if* $\det(\vec{a}', \vec{a}'', \vec{a}''') = 0$.

**Proof:** From the Frenet formulae in (4) and (5) we have

$$\det(\vec{a}', \vec{a}'', \vec{a}''') = \begin{cases} -k_2^5 \left(\dfrac{k_1}{k_2}\right)'; & N \text{ is timelike,} \\ k_2^5 \left(\dfrac{k_1}{k_2}\right)'; & N \text{ is spacelike.} \end{cases} \quad (14)$$

Let now $N$ be a $q$-slant ruled surface in $E_1^3$. By Theorem 3.1, we have $k_1/k_2$ is constant. Then from (14) it follows that $\det(\vec{a}', \vec{a}'', \vec{a}''') = 0$.

Conversely, if $\det(\vec{a}', \vec{a}'', \vec{a}''') = 0$, since the curvature $k_2$ is non-zero from (14) it is obtained that $k_1/k_2$ is constant and Theorem 3.1 gives that $N$ is a $q$-slant ruled surface in $E_1^3$.

***Theorem 3.4.*** *Non-null ruled surface $N$ is a $q$-slant ruled surface if and only if*

$$\vec{q}''' + m\vec{q}' = 3k_1' \vec{h}', \quad (15)$$

*holds where*

$$m = \begin{cases} -\dfrac{k_1''}{k_1} + \varepsilon(k_1^2 - k_2^2); & N \text{ is timelike,} \\ -\left(\dfrac{k_1''}{k_1} + k_1^2 + k_2^2\right); & N \text{ is spacelike.} \end{cases}$$

**Proof:** Assume that $N$ is a timelike $q$-slant ruled surface. From (4) we get

$$\vec{q}'' = -\varepsilon k_1^2 \vec{q} + k_1' \vec{h} + k_1 k_2 \vec{a}, \quad (16)$$

$$\vec{q}''' = (-3\varepsilon k_1 k_1')\vec{q} + (k_1'' + \varepsilon k_1 k_2^2)\vec{h} + (2k_1' k_2 + k_1 k_2')\vec{a} - (\varepsilon k_1^2)\vec{q}'. \tag{17}$$

Since $N$ is a $q$-slant ruled surface, $k_1/k_2$ is constant and by differentiation we have

$$k_1 k_2' = k_2 k_1', \tag{18}$$

and from (4)

$$\vec{h} = \frac{1}{k_1}\vec{q}'. \tag{19}$$

Substituting (18) and (19) in (17) gives

$$\vec{q}''' = \left(\frac{k_1''}{k_1} + \varepsilon(k_2^2 - k_1^2)\right)\vec{q}' + 3k_1'(-\varepsilon k_1 \vec{q} + k_2 \vec{a}). \tag{20}$$

Using the second equation of (4), (15) is obtained from (20).

Conversely, let us assume that (15) holds. Differentiating (19) we obtain

$$\vec{h}' = -\left(\frac{k_1'}{k_1^2}\right)\vec{q}' + \left(\frac{1}{k_1}\right)\vec{q}'', \tag{21}$$

and so,

$$\vec{h}'' = -\left(\frac{k_1'}{k_1^2}\right)'\vec{q}' - 2\left(\frac{k_1'}{k_1^2}\right)\vec{q}'' + \left(\frac{1}{k_1}\right)\vec{q}'''. \tag{22}$$

Substituting (15) in (22) it follows

$$\vec{h}'' = -2\left(\frac{k_1'}{k_1^2}\right)\vec{q}'' - \left[\left(\frac{k_1'}{k_1^2}\right)' + \frac{m}{k_1}\right]\vec{q}' + 3\left(\frac{k_1'}{k_1}\right)\vec{h}'. \tag{23}$$

Now, writing (16) in (23) and using (4) we have

$$\vec{h}'' = -\left[\left(\frac{k_1'}{k_1^2}\right)' + \frac{m}{k_1}\right]\vec{q}' - (\varepsilon k')\vec{q} - 2\left(\frac{k_1'}{k_1}\right)^2 \vec{h} + \left(\frac{k_2 k_1'}{k_1}\right)\vec{a}. \tag{24}$$

On the other hand, from (4) it is obtained

$$\vec{h}'' = -(\varepsilon k_1)\vec{q}' - (\varepsilon k_1')\vec{q} + (\varepsilon k_2^2)\vec{h} + k_2'\vec{a}. \tag{25}$$

Substituting (25) in (24) we have

$$\frac{k_2'}{k_2} = \frac{k_1'}{k_1}. \tag{26}$$

Integrating (26) we get that $k_1/k_2$ is constant and by Theorem 3.1, $N$ is a $q$-slant ruled surface.

If $N$ is a spacelike ruled surface, then by the similar way it is obtained that $N$ is a $q$-slant ruled surface if and only if (15) holds for $m = -\left(\frac{k_1''}{k_1} + k_1^2 + k_2^2\right)$.

**Theorem 3.5.** *Let $N$ be a developable non-null ruled surface in $E_1^3$. Then $N$ is a $q$-slant ruled surface if and only if the striction line $\vec{c}(s)$ is a general helix in $E_1^3$.*

**Proof:** Since $N$ is a developable non-null ruled surface in $E_1^3$, from Theorem 2.1 we have $\vec{c}'(s) = \vec{t}(s) = \vec{q}(s)$ where $\vec{t}(s)$ is the unit tangent of $\vec{c}(s)$. Then from Definition 3.1, it is clear that $N$ is a $q$-slant ruled surface if and only if the striction line $\vec{c}(s)$ is a general helix in $E_1^3$.

## 4. $h$-Slant Ruled Surfaces in $E_1^3$

In this section, we introduce the definition and characterizations of $h$-slant non-null ruled surfaces in $E_1^3$. First, we give the following definition.

**Definition 4.1.** Let $N$ be a non-null ruled surface in $E_1^3$ given by the parametrization
$$\vec{r}(s,v) = \vec{c}(s) + v\vec{q}(s), \quad \|\vec{q}(s)\| = 1, \tag{27}$$
where $\vec{c}(s)$ is striction curve of $N$ and $s$ is arc length parameter of $\vec{c}(s)$. Let the Frenet frame and non-zero invariants of $N$ be $\{\vec{q}, \vec{h}, \vec{a}\}$ and $k_1, k_2$, respectively. Then, $N$ is called a $h$-slant ruled surface if the central normal vector $\vec{h}$ makes a constant angle with a fixed non-zero unit direction $\vec{u}$ in the space, i.e.,
$$\langle \vec{h}, \vec{u} \rangle = c_h = constant. \tag{28}$$

Then, under the assumptions given in Definition 4.1, we can give the following theorems characterizing $h$-slant ruled surfaces.

***Theorem 4.1.*** *$N$ is a non-null $h$-slant ruled surface if and only if the function*
$$f = \begin{cases} \dfrac{k_1^2}{(\varepsilon(k_2^2 - k_1^2))^{\frac{3}{2}}} \left(\dfrac{k_2}{k_1}\right)'; & N \text{ is timelike}, \\[2ex] \dfrac{k_1^2}{(k_1^2 + k_2^2)^{\frac{3}{2}}} \left(\dfrac{k_2}{k_1}\right)'; & N \text{ is spacelike}. \end{cases} \tag{29}$$
*is constant.*

**Proof:** Assume that $N$ is a non-null $h$-slant ruled surface and let $N$ be timelike. Let $\vec{u}$ be a fixed constant vector such that $\langle \vec{h}, \vec{u} \rangle = c_h = constant$. Then for the vector $\vec{u}$ we have
$$\vec{u} = b_1(s)\vec{q}(s) + c_h \vec{h}(s) + b_2(s)\vec{a}(s), \tag{30}$$
where $b_1 = b_1(s)$ and $b_2 = b_2(s)$ are smooth functions of arc length parameter $s$. Since $\vec{u}$ is constant, differentiation of (30) gives
$$\begin{cases} b_1' - \varepsilon c_h k_1 = 0, \\ b_1 k_1 + \varepsilon b_2 k_2 = 0, \\ b_2' + c_h k_2 = 0. \end{cases} \tag{31}$$
From the second equation of system (31) we have
$$b_1 = -\varepsilon b_2 \frac{k_2}{k_1}. \tag{32}$$
Moreover,
$$\langle \vec{u}, \vec{u} \rangle = \varepsilon b_1^2 + c_h^2 - \varepsilon b_2^2 = constant. \tag{33}$$
Substituting (32) in (33) gives
$$\varepsilon b_2^2 \left( \left(\frac{k_2}{k_1}\right)^2 - 1 \right) = n^2 = constant. \tag{34}$$

Then from (34) it is obtained that

$$b_2 = \pm \frac{n}{\sqrt{\left|\varepsilon\left[\left(\frac{k_2}{k_1}\right)^2 - 1\right]\right|}}. \tag{35}$$

Considering the third equation of system (31), from (35) we have

$$\frac{d}{ds}\left[\pm \frac{n}{\sqrt{\left|\varepsilon\left[\left(\frac{k_2}{k_1}\right)^2 - 1\right]\right|}}\right] = -c_h k_2.$$

This can be written as

$$\frac{k_1^2}{\left(\varepsilon(k_2^2 - k_1^2)\right)^{\frac{3}{2}}}\left(\frac{k_2}{k_1}\right)' = \frac{c_h}{n} = constant,$$

which is desired.

Conversely, assume that $N$ is timelike and the function in (29) is constant, i.e.,

$$\frac{k_1^2}{\left(\varepsilon(k_2^2 - k_1^2)\right)^{\frac{3}{2}}}\left(\frac{k_2}{k_1}\right)' = constant = d.$$

We define

$$\vec{u} = \frac{k_2}{\sqrt{|\varepsilon(k_2^2 - k_1^2)|}}\vec{q} - d\vec{h} - \frac{\varepsilon k_1}{\sqrt{|\varepsilon(k_2^2 - k_1^2)|}}\vec{a}. \tag{36}$$

Differentiating (36) with respect to $s$ and using (29) we have $\vec{u}' = 0$, i.e., $\vec{u}$ is a constant vector. On the other hand $\langle \vec{h}, \vec{u} \rangle = constant$. Thus, $N$ is a $h$-slant ruled surface in $E_1^3$.

If $N$ is considered as a spacelike ruled surface, then making the similar calculations, it is obtained that $N$ is a $h$-slant ruled surface if and only if the function $\dfrac{k_1^2}{\left(k_1^2 + k_2^2\right)^{\frac{3}{2}}}\left(\dfrac{k_2}{k_1}\right)'$ is constant.

At the following theorem we give a special case for which the first curvature $k_1$ is equal to 1 and obtain the second curvature for the ruled surface $N$ to be a $h$-slant ruled surface.

***Theorem 4.2.*** *Let $N$ be a non-null in $E_1^3$ with first curvature $k_1 \equiv 1$. Then the central normal vector $\vec{h}$ makes a constant angle $\theta$ with a fixed non-null direction $\vec{u}$, i.e., $N$ is a $h$-slant ruled surface if and only if the second curvature is given as follows*

*i) If $N$ is a timelike ruled surface then*

$$k_2(s) = \pm \frac{s}{\sqrt{|s^2 + \varepsilon \mu^2(\theta)|}}, \tag{37}$$

*where*

$$\mu(\theta) = \begin{cases} \coth\theta; & \text{if } \vec{u} \text{ is a timelike vector,} \\ \tanh\theta; & \text{if } \vec{u} \text{ is a spacelike vector.} \end{cases}$$

**ii)** *If $N$ is a spacelike ruled surface then*

$$k_2(s) = \pm \frac{s}{\sqrt{\left|\eta^2(\theta) - s^2\right|}}, \tag{38}$$

*where*

$$\eta(\theta) = \begin{cases} \tanh\theta; & \text{if } \vec{u} \text{ is a timelike vector,} \\ \coth\theta; & \text{if } \vec{u} \text{ is a spacelike vector.} \end{cases}$$

**Proof:** Let $N$ be a timelike ruled surface with first curvature $k_1 \equiv 1$ and let $N$ be a $h$-slant ruled surface. Then for a fixed constant timelike unit vector $\vec{u}$ we have

$$\langle \vec{h}, \vec{u} \rangle = \sinh\theta = constant. \tag{39}$$

Differentiating (39) with respect to $s$ gives

$$\langle -\varepsilon\vec{q} + k_2\vec{a}, \vec{u} \rangle = 0, \tag{40}$$

and from (40) we have

$$\langle \vec{q}, \vec{u} \rangle = \varepsilon k_2 \langle \vec{a}, \vec{u} \rangle. \tag{41}$$

If we put $\langle \vec{a}, \vec{u} \rangle = \varepsilon x$, we can write

$$\vec{u} = (\varepsilon k_2 x)\vec{q} + (\sinh\theta)\vec{h} - x\vec{a}. \tag{42}$$

Since $\vec{u}$ is unit, from (42) we have

$$x = \pm \frac{\cosh\theta}{\sqrt{\left|\varepsilon(1 - k_2^2)\right|}}. \tag{43}$$

Then, the vector $\vec{u}$ is given as follows

$$\vec{u} = \pm \frac{\varepsilon k_2 \cosh\theta}{\sqrt{\left|\varepsilon(1 - k_2^2)\right|}} \vec{q} + (\sinh\theta)\vec{h} \mp \frac{\cosh\theta}{\sqrt{\left|\varepsilon(1 - k_2^2)\right|}} \vec{a}. \tag{44}$$

Differentiating (40) with respect to $s$, it follows

$$\langle -\varepsilon(1 - k_2^2)\vec{h} + k_2'\vec{a}, \vec{u} \rangle = 0. \tag{45}$$

Writing $\langle \vec{a}, \vec{u} \rangle = \varepsilon x$ and (39) in (45) we have

$$x = \frac{(1 - k_2^2)\sinh\theta}{k_2'}. \tag{46}$$

From (43) and (46) we obtain the following differential equation,

$$\pm \coth\theta \frac{\varepsilon k_2'}{\left(\varepsilon(1 - k_2^2)\right)^{3/2}} - 1 = 0. \tag{47}$$

By integration from (47) we get

$$\pm \coth\theta \frac{\varepsilon k_2}{\sqrt{\left|\varepsilon(1 - k_2^2)\right|}} - s + c = 0, \tag{48}$$

where $c$ is integration constant. The integration constant can be subsumed thanks to a parameter change $s \to -s + c$. Then (48) can be written as

$$\pm\coth\theta\frac{\varepsilon k_2}{\sqrt{\left|\varepsilon(1-k_2^2)\right|}}=-s \qquad (49)$$

which gives us $k_2(s)=\pm\dfrac{s}{\sqrt{\left|s^2+\varepsilon\coth^2\theta\right|}}$.

Conversely, assume that $k_2(s)=\pm\dfrac{s}{\sqrt{\left|s^2+\varepsilon\coth^2\theta\right|}}$ holds and let us put

$$x=\pm\frac{\cosh\theta}{\sqrt{\left|\varepsilon(1-k_2^2)\right|}}=\pm\frac{\cosh\theta}{\sqrt{\left|\varepsilon-\dfrac{\varepsilon s^2}{s^2+\varepsilon\coth^2\theta}\right|}}=\pm\sinh\theta\sqrt{\left|s^2+\varepsilon\coth^2\theta\right|}, \qquad (50)$$

where $\theta$ is constant. Then, $k_2 x = s\sinh\theta$. Let now consider the vector $\vec{u}$ defined by

$$\vec{u}=\sinh\theta\left(\varepsilon s\vec{q}+\vec{h}\mp\left(\sqrt{\left|s^2+\varepsilon\coth^2\theta\right|}\right)\vec{a}\right) \qquad (51)$$

We will prove that $\vec{u}$ is constant and makes a constant angle $\theta$ with $\vec{h}$. By differentiating (51) and using Frenet formulae we have $\vec{u}'=0$, i.e., the direction of $\vec{u}$ is constant and $\langle\vec{h},\vec{u}\rangle=\sinh\theta=constant$. Then $N$ is a $h$-slant ruled surface.

If we assume that $\vec{u}$ is spacelike then we have $\langle\vec{h},\vec{u}\rangle=\cosh\theta=constant$ and making the similar calculations we obtain $k_2(s)=\pm\dfrac{s}{\sqrt{\left|s^2+\varepsilon\tanh^2\theta\right|}}$. Then we can write (37).

If the ruled surface $N$ is a spacelike ruled surface then following the same procedure it is easily obtained that $N$ is a $h$-slant ruled surface if and only if the second curvature is given by $k_2(s)=\pm\dfrac{s}{\sqrt{\left|\eta^2(\theta)-s^2\right|}}$ where $\eta(\theta)=\tanh\theta$, if $\vec{u}$ is a timelike vector; and $\eta(\theta)=\coth\theta$, if $\vec{u}$ is a spacelike vector.

On the other hand, if the striction line $\vec{c}(s)$ is a geodesic on $N$, then the principal normal vector $\vec{n}$ of $\vec{c}(s)$ and the central normal vector $\vec{h}$ of $N$ coincide. Then, we have the following corollary.

***Corollary 4.1.*** *Let the striction line $\vec{c}(s)$ be a geodesic on $N$. Then $N$ is a non-null $h$-slant ruled surface if and only if the striction line is a slant helix in $E_1^3$.*

If the non-null ruled surface $N$ is developable, then by Theorem 2.1, the Frenet frame $\{\vec{t},\vec{n},\vec{b}\}$ of the striction line $\vec{c}(s)$ coincides with the frame $\{\vec{q},\vec{h},\vec{a}\}$ and we can give the following corollary.

***Corollary 4.2.*** *Let $N$ be a non-null developable surface in $E_1^3$. Then $N$ is a $h$-slant ruled surface if and only if the striction line is a slant helix in $E_1^3$.*

# 5. $a$-Slant Ruled Surfaces in $E_1^3$

In this section we introduce the definition of $a$-slant ruled surfaces in $E_1^3$.

**Definition 5.1.** Let $N$ be a non-null ruled surface in $E_1^3$ given by the parametrization
$$\vec{r}(s,v) = \vec{c}(s) + v\vec{q}(s), \ \|\vec{q}(s)\| = 1,$$
where $\vec{c}(s)$ is striction curve of $N$ and $s$ is arc length parameter of $\vec{c}(s)$. Let the Frenet frame and non-zero invariants of $N$ be $\{\vec{q}, \vec{h}, \vec{a}\}$ and $k_1, k_2$, respectively. Then, $N$ is called a $a$-slant ruled surface if the central tangent vector $\vec{a}$ makes a constant angle with a fixed non-zero direction $\vec{u}$ in the space, i.e.,
$$\langle \vec{a}, \vec{u} \rangle = c_a = constant.$$

From (10) it is clear that a non-null ruled surface $N$ is $a$-slant ruled surface if and only if it is a $q$-slant ruled surface. So, all the theorems given in Section 3 also characterize the $a$-slant ruled surfaces.

After these definitions and characterizations of non-null slant ruled surfaces we can give the followings:

Let $N_1$ and $N_2$ be two non-null ruled surfaces in $E_1^3$ with Frenet frames $\{\vec{q}_1, \vec{h}_1, \vec{a}_1\}$ and $\{\vec{q}_2, \vec{h}_2, \vec{a}_2\}$, respectively. If $N_1$ and $N_2$ have common central normals i.e., $\vec{h}_1 = \vec{h}_2$ at the corresponding points of their striction lines, then $N_1$ and $N_2$ are called Bertrand offsets [9]. Similarly, if $\vec{a}_1 = \vec{h}_2$ at the corresponding points of their striction lines, then the surface $N_2$ is called a Mannheim offset of $N_1$ and the ruled surfaces $N_1$ and $N_2$ are called Mannheim offsets [17]. Considering these definitions we come to the following corollaries:

**Corollary 5.1.** *Let $N_1$ be a $h$-slant ruled surface. Then the Bertrand offsets of $N_1$ form a family of $h$-slant ruled surfaces.*

**Corollary 5.2.** *Let $N_1$ and $N_2$ form a Mannheim offset. Then $N_1$ is a $q$-slant (or $a$-slant) ruled surface if and only if $N_2$ is a $h$-slant ruled surface.*

## References

[1] Ali, A.T., Lopez, R., Slant Helices in Minkowski Space $E_1^3$, J. Korean Math. Soc. 48(1) (2011) 159-167.
[2] Ali, A.T., Turgut, M., Position vector of a time-like slant helix in Minkowski 3-space, J. Math. Anal. Appl. 365 (2010) 559–569.
[3] Ali, A.T., Position vectors of slant helices in Euclidean Space $E^3$, arXiv: 0907.0750v1, [math.DG], 4 Jul 2009.
[4] Barros, M., General helices and a theorem of Lancret, Proc. Amer. Math. Soc. 125(5) (1997) 1503-1509.
[5] Beem, J.K., Ehrlich, P.E., Global Lorentzian Geometry, Marcel Dekker, New York, (1981).
[6] Ekmekçi, N., Hacısalihoğlu, H.H., On Helices of a Lorentzian Manifold, Commun. Fac. Sci. Univ. Ank. Series A1, 45 (1996) 45-50.

[7] Ferrandez, A., Gimenez, A., Lucas, P., Null helices in Lorentzian space forms, Internat. J. Modern Phys. A. 16(30) (2001) 4845-4863.

[8] Izumiya, S., Takeuchi, N., New special curves and developable surfaces, Turk J. Math. 28, (2004) 153-163.

[9] Kasap, E., Kuruoğlu, N., The Bertrand offsets of ruled surfaces in $\mathbb{R}_1^3$, ACTA MATHEMATICA VIETNAMICA, 31(1) (2006) 39-48.

[10] Kocayiğit, H., Önder, M., Timelike curves of constant slope in Minkowski space $E_1^4$, J. Science Techn. Beykent Univ. 1 (2007) 311-318.

[11] Kula, L., Yaylı, Y., On slant helix and its spherical indicatrix, Applied Mathematics and Computation. 169 (2005) 600-607.

[12] Küçük, A., On the developable timelike trajectory ruled surfaces in Lorentz 3-space $IR_1^3$, App. Math. and Comp., 157 (2004) 483-489.

[13] O'Neill, B., Semi-Riemannian Geometry with Applications to Relativity, Academic Press, London, (1983).

[14] Önder, M., Similar Ruled Surfaces with Variable Transformations in Minkowski 3-space, TWMS J. App. Eng. Math. 5(2) (2015) 219-230.

[15] Önder, M., Kocayiğit, H., Kazaz, M., Spacelike helices in Minkowski 4-space $E_1^4$, Ann Univ Ferrara, 56 (2010) 335-343.

[16] Önder, M., Uğurlu, H.H., Frenet Frames and Invariants of Timelike Ruled Surfaces, Ain Shams Eng J., 4 (2013) 507-513.

[17] Önder, M., Uğurlu, H.H., On the Developable Mannheim Offsets of Timelike Ruled Surfaces, Proc. Natl. Acad. Sci., India, Sect. A Phys. Sci. 84(4) (2014) 541-54 8.

[18] Önder, M., Uğurlu, H.H., Frenet Frames and Frenet Invariants of Spacelike Ruled Surfaces, Dokuz Eylul University-Faculty of Engineering Journal of Science and Engineering, 19(57) (2017) 712-722.

[19] Önder, M., Kaya, O., Slant null scrolls in Minkowski 3-space $E_1^3$, Kuwait J. Sci. 43(2) (2016) 31-47.

[20] Struik, D.J., Lectures on Classical Differential Geometry, 2nd ed. Addison Wesley, Dover, (1988).

[21] Uğurlu, H.H., Çalışkan, A., Darboux Ani Dönme Vektörleri ile Spacelike ve Timelike Yüzeyler Geometrisi, Celal Bayar Üniversitesi Yayınları, Yayın No: 0006, (2012).